\documentclass[a4]{article}
\usepackage{amsmath}
\usepackage{amsthm}
\usepackage[psamsfonts,mathscr]{euscript}
\usepackage{amssymb}

\begin{document}
\title{Strong Convergence on Weakly Logarithmic Combinatorial
 Assemblies}

\author{E. Manstavi\v cius \thanks{
{ Vilnius University; Address: Naugarduko str. 24, LT-03225 Vilnius,
Lithuania} }}
 \maketitle

\footnotetext{{\it AMS} 2000 {\it subject classification.} Primary
60C05;      secondary 05A16, 60F15.}





\begin{abstract}
We deal with the random combinatorial structures called assemblies.
 By weakening the logarithmic condition which  assures regularity of the number of
 components of a given order, we extend the notion of logarithmic assemblies.
 Using the author's analytic approach, we generalize the so-called
Fundamental Lemma giving independent process approximation in the
total variation distance of the component structure of an assembly.
To evaluate the influence of strongly dependent large components, we
obtain estimates of the appropriate conditional probabilities by
unconditioned ones. These estimates are applied to examine additive
functions defined on such a class of structures. Some analogs of
Major's and Feller's theorems which concern almost sure behavior of
sums of independent random variables are proved.
\end{abstract}




\newtheorem*{cor1}{Corollary 1}


\newtheorem*{cor3}{Corollary 3}

\newtheorem*{FLemma}{Theorem (Fundamental Lemma)}
\newtheorem*{ex}{Example}
\newtheorem{cor}{Corollary}
\newtheorem{thm}{Theorem}
\newtheorem{lem}{Lemma}
\newtheorem{prop}{Proposition}

\def\E{\mathbf{E}}
\def\C{\mathbb{C}}
\def\D{\mathbb{D}}
\def\G{\mathbb{G}}
\def\V{\mathbf{V}}
\def\Ra{\Rightarrow}
\def\N{\mathbb{N}}
\def\R{\mathbb{R}}
\def\S{\mathbb{S}_n}
\def\Z{\mathbb{Z}}
\def\k{\kappa}
\def\e{\varepsilon}
\def\n{$n\to\infty$}
\def\cF{\mathcal F}
\def\cL{\mathcal L}
\def\cK{\mathcal K}

\def\cA{{\mathcal A}}
\def\re{{\rm e}}
\def\rO{{\rm O}}
\def\ro{{\rm o}}
\def\rd{{\rm d}}

\def\s{\smallskip}
\def\b{\bigskip}


\section{Introduction}
\label{Intr}

In part, this work was stimulated by a critical remark made by R.
Arratia, A.D. Barbour and S. Tavar\'e \cite{ABT1} about analytic
methods applied in the theory of random combinatorial structures. On
page 1622  they wrote: {\it In contrast} (to their method), {\it the
complex analytic approaches typically require conditions to be
satisfied that can be verified in the well-known examples, but which
are difficult to express directly in terms of the basic parameters
of the structures}. Such was the  criticism  to the method
cultivated in the papers by P.~Flajolet and M.~Soria \cite{F-S} and
J.~Hansen \cite{H}. The works written by D.~Stark \cite{St1} and
\cite{St2} could be added to this list as well. Indeed, the
conditions posed on the generating series of structure classes have
some disadvantages.

The authors of \cite{ABT1} did not notice the broader possibilities
hidden in the analytic approach proposed in our papers \cite{EM-96},
\cite{EM-98}, \cite{EM-02},  and refined in \cite{J-E-V-07} and
\cite{EM-09}. So far, this approach was applied to obtain asymptotic
formulas for some Fourier transforms of distributions. That led to
general one-dimensional limit theorems, including the optimal
remainder term estimates. In this regard, apart from the above
mentioned, the papers by V.~Zacharovas \cite{VZ1}, \cite{VZ2}, and
\cite{VZ3} were noticeable. On the other hand, there exist a lot of
works dealing with the deeper total variation approximation (see,
for instance, \cite{ABT} and the references therein). The main goal
of the present paper is to demonstrate that such total variation
approximations can be
 obtained by our method and, at the same time, under more general conditions possed on
  {\it the basic parameters} of the structures. For simplicity, we confine
  ourselves to classes of {\it
 assemblies} or {\it abelian partitional complexes} (see \cite{F}).
 For completeness, we recall the definition and some
properties which can be found in \cite{ABT}.

Let  $\sigma$ be a set of $n\geq1$ points, partitioned into subsets
so that there are $k_j(\sigma)>0$ subsets of size $j$, $1\leq j\leq
n$ and $\bar k(\sigma):=\big(k_1(\sigma),\dots,k_n(\sigma)\big)$. If
$\ell(\bar s):=1s_1+\cdots+ns_n$, where  $\bar
s=(s_1,\dots,s_n)\in\mathbb{Z}_+^n$, then $\ell\big(\bar
k(\sigma)\big)=n$.
 Assume that in each such subset of size $1\leq j\leq n$ by some rule one of  $0<
 m_j<\infty$
possible structures can be chosen. A subset with a structure is a
{\it component} of $\sigma$, and the set $\sigma$ itself  is called
an assembly \cite{ABT}. Using  all possible partitions of $\sigma$
and the same rule to define a structure in a component, we get the
class $ {\cal A}_n$ of assemblies of size $n$. Let $ {\cal A}_0 $ be
comprised of the empty set. The union
\[
    {\cal A}_0\cup {\cal A}_1\cup\cdots\cup {\cal A}_n\cup\cdots
\]
 forms the whole class of assemblies. Its  basic parameters appear in the conditions posed
on the sequence $m_j$, $j\geq1$.

    There are
$$
n!\prod_{j=1}^n \left({1 \over j!} \right)^{s_j}{1 \over s_j!}
$$
 ways to partition an $n$-set into subsets, so that $\bar k(\sigma)=\bar s$ if $\ell(\bar s)=n$ and
  $\bar
s\in\mathbb{Z}_+^{n}$. Hence, there are
$$
Q_n(\bar s):=n!\prod_{j=1}^n \left({m_j\over j!}
\right)^{s_j}{1\over s_j!}
$$
assemblies with the  component vector $\bar k(\sigma)=\bar s$, and
the total number of them in the class $ {\cal A}_n$ equals
\[
           |{\cal A}_n|=\sum_{\ell(\bar s) =n}Q_n(\bar s).
\]

 On the class $ {\cal A}_n$, one can define
the uniform probability measure denoted by
\[
 \nu_n(\dots)=|{\cal A}_n|^{-1}|\{\sigma\in {\cal A}_n, \dots\}|.
\]
From now $\sigma\in {\cal A}_n$ is an elementary event. Following
the tradition of probabilistic number theory and in contrast to
\cite{ABT}, we prefer to leave it defining random variables (r.vs)
on ${\cal A}_n$. The component vector $\bar k(\sigma)$ has the
following  distribution:
\begin{eqnarray*}
\nu_n(\bar{k}(\sigma)=\bar s)={\bf{1}} \{\ell(\bar s) =n\}{{n!}\over
{|{\cal A}_n|}}\prod_{j=1}^n {1\over {s_j!}}\left({m_j\over
j!}\right)^{s_j},
\end{eqnarray*}
where $\bar s=(s_1,\dots,s_n)\in\Z_+^n$.
 This leads to the {\it Conditioning Relation} (see
\cite{ABT}, page 48)
 \begin{eqnarray}
 \nu_n( \bar{k}(\sigma)=\bar s)=P\big(\bar \xi=\bar s|\ell(\bar \xi)=n),
\label{CR}
 \end{eqnarray}
 where $\bar\xi:= (\xi_1,\dots,\xi_n)$ and
 $\xi_j, j\geq 1$, are mutually
independent Poisson r.vs defined on some probability space
$\{\Omega, \cal{F},P\}$ with $\mathbf{E}\xi_j =u^jm_j/j!$, $j\geq1$,
 where $u>0$ is an arbitrary  number.

The so-called {\it Logarithmic Condition} (see \cite{ABT}) in the
case of assemblies requires that
$$
       m_j/{j!}\sim \theta y^{j}/j
$$
for some constants $ y>0$ and  $\theta>0$ as $j\to\infty$. Under
this condition, it is natural and technically convenient to take
$u=y^{-1}$, which yields the relation  $\mathbf{E}\xi_j\sim\theta/j$
as $j\to\infty$.

  Generalizing the Ewens probability in the symmetric group of permutations,
  the author in \cite{EM-02} and \cite{EM-09} examined  random assemblies
taken with weighted frequencies. The research was extended by V.
Zacharovas \cite{VZ4}. Going along this path, one can take a
positive sequence $w_j$, $j\geq 1$, and define
\[
     w(\sigma)=\prod_{j=1}^n w_j^{k_j(\sigma)}, \qquad
     W_n=\sum_{\sigma\in\cA_n} w(\sigma).
\]
Further, one can introduce the probability measure $\nu_n^{(w)}$ on
$\cA_n$ by
  \[
       \nu_n^{(w)}\big(\{\sigma\}\big)= w(\sigma)/W_n, \quad \sigma\in\cA_n.
       \]
Conditioning Relation  (\ref{CR}) still holds for $\nu_n^{(w)}$
instead of $\nu_n$ with the  poissonian random vector $\bar\xi$
provided that $\E \xi_j=u^jm_jw_j/j!$, where $j\geq 1$ and $u>0$ is
an arbitrary constant.  Having all this in mind,
 we extend  the  logarithmic class of
 assemblies discussed in \cite{ABT} and in many previous papers.

\s

 {\bf Definition.} {\it Let $n\geq 1$ and let $\mu_n$ be
a probability measure on $\cA_n$. The pair $\big(\cA_n, \mu_n\big)$
will be called weakly logarithmic  if there exists a random vector
$\bar\xi= (\xi_1,\dots,\xi_n)$ with mutually independent poissonian
coordinates such that
\begin{eqnarray*}
\mu_n(\bar{k}(\sigma)=\bar{s})=P\big(\bar \xi=\bar s|\ell(\bar
\xi)=n\big)
 \end{eqnarray*}
for each $\bar s\in\Z_+^n$ and
\begin{equation}
                 { \theta' \over j}\leq \lambda_j:=\mathbf{E}\xi_j\leq
                 {\theta''\over j}
\label{Def}
\end{equation}
 uniformly in $j\geq1$ for some positive
constants  $\theta'$ and $\theta''$.} \s

In our notation,  the {\it logarithmic assemblies} are characterized
by the condition $\lambda_j\sim \theta/j$ as $j\to\infty$, where
$\theta>0$ is a constant (see \cite{ABT}).

 The main result
of this paper is the following total variation approximation. Let
$\cL(X)$ be the distribution of a r.v. $X$. Afterwards the index
$r$, $1\le r\le n$, added to the vectors $\bar k(\sigma)$ and $\bar
\xi$ will denote that only the first $r$ coordinates are taken. Let
$x_+=\max\{x,0\}$ for $x\in\R$ and $\ll$ be an analog of the symbol
$\rO(\cdot)$.

\begin{FLemma} Let  $(\cA_n,\mu_n)$ be weakly logarithmic.
There exist positive constants $c_1$ and
 $c_2$ depending on $\theta'$ and such that
\begin{equation}
 \rho_{TV}\Big(\cL\big(\bar k_r(\sigma)\big),
\cL(\bar\xi_r)\Big):= \sum_{\bar s\in\Z_+^r} \Big(\mu_n\big(\bar
     k_r(\sigma)=\bar s)-P(\bar\xi_r=\bar s)\Big)_+\ll \Big({r\over n}\Big)^{c_1}
\label{TV}
\end{equation}
uniformly in $1\leq r\leq c_2 n$. The constant in $\ll$ depends on
$\theta'$ and $\theta''$ only.
\end{FLemma}

Adopting I. Z. Ruzsa's idea going back to probabilistic number
theory (see \cite{Ru}),  we \cite{EM-98LMJ} observed that some
conditional discrete probabilities can be estimated by appropriate
unconditional ones. This led to upper estimates of the distributions
$\cL\big(\bar k(\sigma)\big)$ of the cycle  structure vector $\bar
k(\sigma)$ of a random permutation $\sigma$ under the uniform
probability defined on the symmetric group. In the joint paper with
G.J. Babu \cite{EJ-99}, the idea was extended to permutations taken
with the Ewens probability and later, jointly with J. Nork\=unien\.e
\cite{EM-JN}, we adopted it for logarithmic assemblies. We now
develop the same principle for weakly logarithmic assemblies.

Firstly, we introduce some notation in the semi-lattice $\Z_+^{n}$
taken from the theory of euclidean spaces. For two vectors $\bar
s=(s_1,\dots,s_n)$ and $\bar t=(t_1,\dots,t_n)$, we set $\bar
s\perp\bar t$ if $s_1t_1+\cdots+s_nt_n=0$ and  write $\bar s\leq
\bar t$ if $s_j\leq t_j$ for each $j\leq n$. Further, we adopt the
notation $\bar s\parallel\bar t$ for the expression ``$\bar s$
exactly enters $\bar t$'' which means that $\bar s\leq\bar t$ and
$\bar s\perp \bar t-\bar s$. For arbitrary subset
$U\subset\Z_+^{n}$, we define its extension
\begin{equation}
V=V(U)=\big\{\bar s=\bar{ t^1}+\bar{ t^2}-\bar{ t^3}:\;\bar{
t^1},\bar{t^2},\bar{ t^3}\in U,\, \bar{ t^1}\perp (\bar {t^2}-\bar
{t^3}),\, \bar{t^3}\parallel \bar{t^2}\big\}. \label{VU}
\end{equation}
Set also $\overline A=\Z_+^{n}\setminus A$ and $\theta=\min\{1,
\theta'\}$.

\begin{thm} \label{T-2} Let  $(\cA_n,\mu_n)$ be weakly logarithmic and $\bar \xi$ be
 the poissonian random vector introduced in Definition. For
arbitrary  $U\in \Z_+^{n}$,
\[
    \mu_n\big(\bar k(\sigma)\in \overline V\big)=P\big(\bar \xi\in \overline V\vert\, \ell(\bar\xi)=n\big)
    \ll P^{\theta}(\bar
    \xi\in\overline U)+{\bold 1}\{\theta<1\}n^{-\theta},
\]
where the implicit constants depend on $\theta'$ and $\theta''$
only.
\end{thm}

   The claim of Theorem \ref{T-2} becomes more
   transparent when applied to the value distributions of  additive
   functions. We demonstrate this in a fairly general context.
   Let $(\G,+)$ be an abelian group and $h_j(s)$, $j\in\N$,
   $s\in\Z_+$, be a two-dimensional sequence in $\G$
   satisfying the condition $h_j(0)=0$ for each $j\geq1$. Then
we can define an
    {\it additive function} $h\colon  {\cal A}_n\to\mathbb   G$ by
      \begin{equation}
            h(\sigma)=\sum_{j\le n} h_j\big(k_j(\sigma)\big).
      \label{h}
      \end{equation}
If $h_j(s)=a_j s$ for some $a_j\in \G$, where $j\in \N$ and
$s\in\Z_+$, then the function $h$ is called  {\it completely
additive}.

\begin{cor} \label{cor-1} Let $(\G,+)$ be an abelian group and $h\colon  {\cal A}_n\to\mathbb   G$
be an additive function. Uniformly in  $A\subset \G$,
\[
     \mu_n\big(h(\sigma)\not\in A+A-A\big)
    \ll P^{\theta}\bigg(\sum_{j\leq n} h_j(\xi_j)\not\in A\bigg)+{\bold 1}\{\theta<1\}n^{-\theta}.
\]
\end{cor}

\begin{cor} \label{cor-2} Let  $h\colon  {\cal A}_n\to\mathbb   R$
be an additive function. Uniformly in  $a\in\R$ and $u\geq0$,
\[
     \mu_n\big(|h(\sigma)-a|\geq u\big)
    \ll P^{\theta}\bigg(\bigg|\sum_{j\leq n} h_j(\xi_j)-a\bigg|\geq u/3\bigg)+{\bold
    1}\{\theta<1\}n^{-\theta}.
\]
\end{cor}

As in the case of logarithmic assemblies, Fundamental Lemma and
Theorem 1 can be used to prove general limit theorems for additive
functions defined on $\cA_n$. One can deal with the one-dimensional
case (see, for instance, \cite{ABT}, Section 8.5) or examine the
weak convergence of random combinatorial processes (see
\cite{EJ-99}, \cite{EJ-02}, \cite{EJ-02a},  \cite{EM-02Pal}, and
\cite{ABT}, Section 8.1). This approach can be applied to examine
the strong convergence. Extending papers \cite{EM-04} and
\cite{JN1}, we now obtain an analog of the functional law of
iterated logarithm. It can be compared with Major's \cite{Ma} result
for i.r.vs, generalizing the celebrated Strassen's theorem.

 It is worth stressing that we deal with random variables which are
 defined on a sequence of probability spaces, not on a fixed space.
  This raises the first obstacle
 to be overcome; therefore, we adopt some basic definitions.

  Let $(S,d)$ be a separable metric space.
Assume that $X,X_1, X_2,\dots,X_n$ are $S$-valued random variables
all defined on the probability space $\{\Omega_n, {\cal F}_n,
P_n\}$. Denote by $d(Y,A):=\inf\{d(Y,Z):Z\in A\},\quad A\subset S$,
$Y\in S$, the distance from $Y$ to $A$. We say that $X_m$ converges
to $X$ $\{P_n\}$-almost surely ($\{P_n\}$-a.s.), if for each
$\varepsilon>0$
$$
\lim_{n_1\to\infty}\limsup_{n\rightarrow\infty} P_n\big(\max_{n_1\le
m\leq n}d(X_m,X)\geq \varepsilon \big)=0.
$$
If $P_n=P$ does not depend on $n$, our definition agrees with that
of classical almost sure convergence (see \cite{Pe}, Chapter X).
 A compact set $A\subset S$ is called a {\it cluster}  for the sequence $X_m$ if,  for each $\varepsilon>0$
and each $Y\in A$,
$$
\lim_{n_1\to \infty}\limsup_{n\to \infty}P_n\big(\max_{n_1\leq m\le
n}d(X_m,A)\geq \varepsilon \big)=0
$$
and
$$
\lim_{n_1\to
\infty}\liminf_{n\rightarrow\infty}P_n\big(\min_{n_1\leq m\leq
n}d(X_m,Y)< \varepsilon \big)=1.
$$
We denote the last two relations, by
   $$
    X_m\Rightarrow A, \qquad(\{P_n\}{\textrm{-a.s.}})
   $$

 Let $C[0,1]$ be the Banach space of continuous functions  on the interval
$[0,1]$ with the supremum distance $\rho(\cdot, \cdot)$.
   The set of absolutely continuous functions $g$ such that $g(0)=0$ and
   $$
          \int_0^1(g'(t))^2dt\le1
   $$
is called the {\it Strassen set} $\cK$. We shall show that it is the
cluster set of some combinatorial processes constructed using
partial sums
\[
            h(\sigma,m):=\sum_{j\le m} h_j\big(k_j(\sigma)\big),
      \]
where $ h_j(s)\in \R$ and $ 1\le m\le n$. Set $a_j=h_j(1)$,
\begin{eqnarray*}
  A(m):=\sum_{j=1}^{m}a_j(1-\re^{-\lambda_j}),\qquad
      B^2(m):=\sum_{j=1}^{m}a_j^2e^{-\lambda_j}\big(1-e^{-\lambda_j}\big),
         \end{eqnarray*}
 and $\beta(m)=B(m)\sqrt{2LLB(m)}$, where $1\leq m\leq n$.
We denote by $u_m(\sigma,t)$ the polygonal line joining the points
\[
   (0,0), \qquad \big(B^2(i), h(\sigma,i)-A(i)\big), \quad 1\leq
   i\leq m,
   \]
and set
$$
U_m(\sigma,t)=\beta(m)^{-1}u_m(\sigma,B^2(m)t), \qquad
\sigma\in\cA_n, \; 0\le t\le 1,
$$
for $1\leq m\leq n$.
 The following result generalizes the cases examined in
\cite{EM-98LMJ}, \cite{JN}, and \cite{JN1}.

\begin{thm}\label{T-3} Let $(\cA_n, \mu_n)$ be weakly logarithmic. If $B(n)\to \infty$ and
\begin{equation}
 a_j={\rm
o}\bigg({B(j)\over \sqrt{LL B(j)}}\bigg), \qquad j\to \infty,
\label{c-Kol}
\end{equation}
then
 \begin{eqnarray}
 U_m(\sigma,\cdot)\Rightarrow {\cal K} \qquad (\{\mu_n\}\textrm{-a.s.}).
 \end{eqnarray}
\end{thm}

\smallskip
 Applying continuous functionals defined on the space $C[0,1]$, we derive
 partial cases of the last theorem.

 \begin{cor} Let the  conditions of Theorem $\ref{T-3}$ be satisfied. The following relations hold
 $\{\mu_n\}$-a.s.\,

\begin{itemize}
\item[$(i)$] $U_m(1)\Rightarrow [-1,1];$

\item[$(ii)$]$\big(U_m(\sigma, 1/2), U_m(\sigma, 1)\big)\Rightarrow \{(u,v)\in\R^2:\;u^2+(v-u)^2\le 1/2\};$

\item[$(iii)$] \, if\, $U_{m'}(\sigma, 1/2)\Rightarrow \sqrt{2}/2$ for some
subsequence $m'\to \infty$, then $U_{m'}(\sigma, \cdot)\Rightarrow
g_1$, where
$$
g_1(t)=\left\{%
\begin{array}{ll}
    t\sqrt{2} \quad & {\rm if}\quad 0\le t\le 1/2, \\
    \sqrt{2}/2 \quad & {\rm if}\quad 1/2\le t\le 1; \\
\end{array}%
\right.
$$
\item[$(iv)$] \, if\, $U_{m'}(\sigma, 1/2)\Rightarrow 1/2$ and
$U_{m'}(\sigma, 1)\Rightarrow 0$ for some subsequence $m'\to
\infty$, then $U_{m'}(\sigma, \cdot)\Rightarrow g_2$, where
$$
g_2(t)=\left\{%
\begin{array}{ll}
    t   \quad &{\rm if}\quad 0\le t\le 1/2, \\
    1-t \quad &{\rm if}\quad 1/2\le t\le 1. \\
\end{array}%
\right.
$$
\end{itemize}
\end{cor}

Using other more sophisticated functionals (see, e.g.,  \cite{Fr},
Chapter I),  one can proceed in a similar manner. Claim $(i)$
includes the  assertion that
\[
         \big|h(\sigma,m)-A(m)\big|\le (1+\e)\beta_m
\]
 holds uniformly in $m $, $n_1\leq m\leq n$, for asymptotically almost all
$\sigma\in\mathcal{A}_n$ as $n$ and $n_1$ tend to infinity.
Moreover, it shows that the upper bound is sharp apart from the term
$\e \beta(m)$. An idea how to improve this error goes back to
W.~Feller's paper \cite{Fel}. It has been exploited  by the author
\cite{EM-04} in the case of a special additive function defined on
permutations. Recently, that paper was generalized for the
logarithmic assemblies  \cite{EM-JN}. We now formulate a more
general result.

 We say that an
increasing  sequence $\psi_m$, $m\geq1$, belongs to the {\it upper
class}  $\Psi^+$ (respectively, {\it the lower class} $\Psi^-$) if
\begin{equation}
\lim_{n_1\to\infty}\limsup_{n\to\infty}\mu_n\Big(\max_{n_1\leq m\leq
n} \psi_m^{-1}\big|h(\sigma, m)-A(m)\big|\geq 1\Big)=0,
 \label{uppercl}
\end{equation}
\[
\bigg(\lim_{n_1\to\infty}\liminf_{n\to\infty}\mu_n\Big(\max_{n_1\leq
m\leq n} \psi_m^{-1}\big|h(\sigma, m)-A(m)\big|\geq 1\Big)=1\bigg).
\]

\begin{thm} \label{T-3} Let $(\cA_n, \mu_n)$ be  weakly logarithmic and $B(n)\to\infty$.
 Assume that  a positive sequence $\phi_n\to\infty$ is
   such that
  \begin{equation}
     a_j=\rO\bigg( {B(j)\over \phi_j^3}\bigg), \quad j\ge1.
\label{Kol1}
 \end{equation}
 If the series
\begin{equation}
                  \sum_{j=1}^\infty{a_j^2\phi_j\over
                  jB^2(j)}\,\re^{-\phi_j^2/2}
\label{ser}
\end{equation}
converges, then $B(m)\phi_m \in\Psi^+$. If  series $(\ref{ser})$
diverges, then  $B(m)\phi_m \in\Psi^-$.
\end{thm}

 Since the series
\[
                  \sum_{j=1}^\infty{a_j^2 \over j}{ (LL B(j))^{1/2}\over B^2(j)(L B(j))^{1+x}}
\]
converges for $x=\e$ and diverges for $x=-\e$, the last theorem
implies $(i)$ in Corollary 3 under a bit stronger condition.
 To illustrate Theorem \ref{T-3},
 let     $  \gamma_{2m}^2(\pm\e):=2(1\pm\e)L_2
   B(m)$,
   \[
   \gamma_{3m}^2(\pm\e)/2:=L_2 B(m)+{3\over2}(1\pm\e)L_3 B(m),
   \]
   and
\[
   \gamma_{sm}^2(\pm\e)/2:=L_2 B(m)+{3\over2}L_3 B(m)+L_4 B(m)+\cdots+(1\pm\e)L_s B(m)
   \]
for $s\geq4$.

\begin{cor} Under the conditions of Theorem $\ref{T-3}$,
 we have
\[
         B(m)\gamma_{sm}(\e)\in\Psi^+
\]
          and
\[
         B(m)\gamma_{sm}(-\e)\in\Psi^-
\]
 for each $s\ge2$.
 \end{cor}

 More corollaries, as  in the case of the  logarithmic
assemblies (see \cite{EM-JN}), could be further formulated.
 The main argument  in deriving  Theorems 2 and 3 is the same;
therefore, we will omit  the proofs of the second result and its
corollaries. The technical details in the case of logarithmic
assemblies can be found in \cite{EM-JN}. Finally, we observe that by
substituting r.vs $\xi_j$, $1\leq j\leq n$, by appropriate
independent geometrically distributed and negative binomial r.vs,
one can similarly extend the logarithmic classes of additive
arithmetical semigroups and weighted multisets (see \cite{ABT}).


 \section{Proof of the Fundamental Lemma}
\label{S-2}

The first lemma reduces the problem to a one-dimensional case. For
$\bar s=(s_1,\dots, s_n)$, set $\ell_{ij}(\bar
s)=(i+1)s_{i+1}+\cdots+js_j$ if $0\leq i<j\leq n$. Moreover, let
$\ell_r(\bar s):=\ell_{0r}(\bar s)$, where $1\leq r\leq n$. Then
$\ell_n(\bar s)=\ell(\bar s)$.

\begin{lem}\label{l-1} We have
\begin{eqnarray}
   && \rho_{TV}\Big(\cL\big(\bar k_r(\sigma)\big),\cL(\bar
    \xi_r)\Big)=
\rho_{TV}\Big(\cL \big(\bar \xi_r\big|\ell(\bar \xi)=n\big),
\cL\big(\bar \xi_r\big)\Big)\nonumber\\
&=& \sum_{m\in\Z_+}P\big(\ell_r(\bar
\xi)=m\big)\bigg(1-{P\big(\ell_{rn}(\bar \xi)=n-m\big)\over
P\big(\ell(\bar\xi)=n\big)}\bigg)_+ \label{onedim}
\end{eqnarray}
\end{lem}

{\it Proof} See \cite{ABT}, p. 60.

Consequently, the ratio of probabilities on the right-hand side in
(\ref{onedim}) is now the main objective. So far, the authors
\cite{ABT1}, assuming the Logarithmic Condition, kept obtaining the
limit approximations as $n\to\infty$ for either of the
probabilities, and then showing their equivalence in a fairly large
region for $m$. The limiting behavior of the probabilities can be
rather complicated for weakly logarithmic  assemblies but, as we
will show in the sequel, the ratio of probabilities in
(\ref{onedim}) is regular. Since
\begin{equation}
P\big(\ell_{rn}(\bar \xi)=m\big)={1\over 2\pi i}\int_{|z|=1}{1\over
z^m}\exp\bigg\{\sum_{r<j\leq n}\lambda_j(z^j-1)\bigg\} \rd z,
\label{P}
 \end{equation}
 one can apply  our analytic technique (see
\cite{EM-02} or \cite{EM-09}) which has been elaborated to compare
the Taylor coefficients of two power series. Namely, if $d_j\in\R_+$
and $f_j\in\C$, $1\leq j\leq n$, are two sequences, the latter maybe
depending on $n$ or other parameters, and
  \[
         D(z):=\exp\bigg\{\sum_{j\le n}{d_j\over j}z^j\bigg\}=:\sum_{s=0}^\infty
         D_sz^s,
\]
\[
         F(z):=\exp\bigg\{\sum_{j\le n}{f_j\over j}z^j\bigg\}=:\sum_{s=0}^\infty
         F_sz^s,
\]
then, under certain  conditions, we have obtained asymptotic
formulas for $F_n/D_n$ as $n\to\infty$. As in \cite{EM-02}, we now
also assume the inequalities
\begin{equation}
             d'\leq d_j\leq d''
\label{d}
\end{equation}
 for all $1\leq j\leq n$ and some positive constants $d'\leq d''$.
In our case, $f_j$ are very special; therefore, we can simplify the
previous argument  and get rid of (2.4) in \cite{EM-02}. The goal
now is to find the ratio $F_m/D_n$ preserving some uniformity.

Set, for brevity,
 \[
            e_r=\exp\bigg\{-\sum_{j\leq r}{d_j\over j}\bigg\}.
 \]

\begin{prop}\label{prop1} Assume that the sequence $d_j$, $1\leq j\leq n$, satisfies condition~$(\ref{d})$. For  $0\leq r\leq n$, set $f_j=d_j$ if $r<j\leq n$ and
$f_j=0$ if $j\leq r$. Let $0\leq \eta \leq 1/2$ and $1/n\leq \delta
\leq 1/2$ be arbitrary. There exists  a  positive constant $c$
depending on $d'$ only such that
\[
      F_m/(e_rD_n)-1\ll \big(\eta+(r/n){\bf 1}\{r\geq 1\}\big) \delta^{-1}+\delta^c
\]
uniformly in
\begin{equation}
     0\leq r\leq \delta n, \qquad n(1-\eta)\leq m\leq n.
\label{reg}
\end{equation}
Here and in the proof of this claim, the constant in  $\ll$ depends
on $d'$ and $d''$ only.
\end{prop}

We will use the following notation. Let $K$, $1\leq\delta n<K\leq
n$, be a parameter to be chosen later. For a fixed $0<\alpha<1$, we
introduce the  functions
 \[
 q(z):=\sum_{r<j\leq n}d_jz^{j-1},\qquad
G_1(z)= \exp\bigg\{\alpha\sum_{r<j\le K}{d_j\over j}z^j\bigg\},
\]
\[
G_2(z)= \exp\bigg\{-\alpha\sum_{K <j\leq n}{d_j\over j}z^j\bigg\},
\qquad G_3(z)=F^\alpha(z)-G_1(z).
\]
We denote by $[z^k]U(z)$ the $k$th Taylor
 coefficient of an analytic at zero  function
$U(z)$. Observe that
 \begin{equation}
   [z^k] G_3(z)\leq [z^k] F^\alpha(z),\quad k\geq0,
\label{equal}
\end{equation}
 where $a_j=d_j$ if $r< j\leq n$, and $a_j=0$
otherwise.  Set further $T= (\delta n)^{-1}$,
 \[
 \Delta=\{z=e^{it}:\; T< |t|\leq\pi\}, \qquad
\Delta_0=\{z=e^{it}:\; |t|\leq T\}.
\]

  Seeking $F_m$, we start from the following
identity
\begin{eqnarray}
 F_m&=&{1\over 2\pi i m}\int_{|z|=1}{F'(z)\over z^m}\,\rd
 z\nonumber\\
 &=&{1\over 2\pi i m}\bigg(\int_{\Delta_0}+\int_{\Delta}\bigg){F'(z)\big(1-G_2(z)\big)\over z^m}\,
 \rd z\nonumber\\
 &&\quad +{1\over 2\pi i m}\int_{|z|=1}{F'(z)G_2(z)\over z^m}\,\rd z
 =: J_0+J_1+J_2. \label{J012}
\end{eqnarray}
In what follows, we estimate the integrals $J_1$ and $J_2$ and,
changing the integrand, reduce $J_0$ to the main term of an
asymptotical formula for $D_n$. The proof of Proposition \ref{prop1}
consists of a few lemmas.

\begin{lem} \label{l-2} We have
\[
    D(1)n^{-1}
\ll D_n\ll D(1)n^{-1}
\]
for all $n\ge1$.
\end{lem}

{\it Proof.} This is Lemma 3.1 from \cite{EM-02}.

\s

\begin{lem} \label{l-3} If $0<\alpha<1$ and $\delta n\geq1$, then
\[
    J_2\ll D_ne_r(K/n)^{\alpha d'}
\]
uniformly in $n/2\leq m\leq n$.
\end{lem}

{\it Proof}. For brevity, let
\[
 u_s:=[z^s]G_1(z), \qquad v_l:=[z^l]F^{1-\alpha}(z), \quad s, l\geq0.
\]
Since
\[
   F'(z)G_2(z)=q(z)G_1(z)F^{1-\alpha}(z),
   \]
   from Cauchy's formula, we have
\[
 J_2={1\over 2\pi i m}\int_{|z|=1}q(z)G_1(z)F^{1-\alpha}(z) {\rd
   z\over z^m}
   = {1\over m}\sum_{r<j\leq m}d_j\sum_{s+l=m-j} u_s v_l.
   \]
Hence, by condition (\ref{d}),
 \begin{eqnarray*}
   J_2&\leq& {2d''\over n}\sum_{s\leq n}u_s  \sum_{l\leq n}v_l\\
   &\leq& {2d''\over n}\,F^{1-\alpha}(1) G_1(1)={2d''F(1)\over n}
   \exp\bigg\{-\alpha\sum_{K< j\leq n}{d_j\over j}\bigg\}\\
   &\ll& D_n e_r(K/n)^{\alpha d'}.
\end{eqnarray*}
In the last step we used Lemma \ref{l-2}.

  The lemma is proved.

\begin{lem} \label{l-4} Let $\delta n\geq 1$. Then
\[
               \max_{T\leq |t|\leq \pi}|F(\re^{it})| \ll
               e_rD(1)\delta^{d'}
               \]
uniformly in $0\leq r\leq \delta n$.
\end{lem}

{\it Proof}. By definition,
\begin{eqnarray}
   {|F(\re^{it})|\over D(1)}&=& e_r{|F(\re^{it})|\over F(1)}=
e_r\exp\bigg\{\sum_{r<j\leq n} {d_j(\cos t j-1)\over j}\bigg\}\nonumber\\
&\leq& e_r\exp\bigg\{d'\sum_{\delta n<j\leq n} {\cos t j-1\over
j}\bigg\} \label{F/D}
\end{eqnarray}
 uniformly in $0\leq r\leq \delta n$. We now use the
relation
\[
S(x, t):=\sum_{j\le x}{\cos tj-1\over j}=\log \min \Big\{1,
{2\pi\over x|t|}\Big\}+\rO(1),
\]
valid for all $x\geq 1$ and $|t|\leq\pi$. It shows that $S(\delta n,
t)\ll 1$  for $T=(\delta n)^{-1}\leq |t|\leq \pi$. Hence, for such
$t$,
\[
   S(n,t)-S(\delta n, t)\leq S(n,T)+\rO(1)=\log \delta+\rO(1).
\]
 This yields the desired claim.

\begin{lem} \label{l-5} Let  $0<\alpha<1$ be arbitrary and $\delta n\geq1$. Then
\[
           J_1\ll {e_rnD_n\over K}\delta^{d'(1-\alpha)}
\]
uniformly in $n/2\leq m \leq n$ and $0\leq r\leq \delta n$.
\end{lem}

{\it Proof}. Recalling the previous notation, we can rewrite
\[
  J_1 ={1\over 2\pi i m}\int_{\Delta}q(z)F^{1-\alpha}(z)G_3(z)  {{\rd}z\over z^{m}}.
  \]
  Hence, by Lemma \ref{l-4},
\begin{eqnarray*}
  J_1 &\ll& n^{-1}\max_{z\in \Delta}|F(z)|^{1-\alpha}
  \int_{|z|=1}\big|q(z)\big|\big|G_3(z)\big|
  |{\rd}z|\\
  &\ll&n^{-1}\Big(e_rD(1)\delta^{d'}\Big)^{1-\alpha}\bigg(\int_{|z|=1}\big|q(z)\big|^2
  |{\rd}z|\bigg)^{1/2}\\
  &&\quad \times\bigg(\int_{|z|=1}|G_3(z)|^2
  |{\rd}z|\bigg)^{1/2}.
\end{eqnarray*}
By Parseval's equality,
\[
    \int_{|z|=1}\big|q(z)\big|^2
  |{\rd}z|=2\pi\sum_{r< j\leq n} d_j^2\leq 2\pi (d'')^2n
  \]
and, recalling (\ref{equal}),
\begin{eqnarray*}
&&\int_{|z|=1}|G_3(z)|^2
  |{\rd}z|\leq 2\pi \sum_{l>K}\big([z^l]G_3(z)\big)^2\\
  &\leq& {2\pi\over
  K^2}\sum_{l=1}^\infty l^2 \big([z^l]F^\alpha(z)\big)^2\ll
      {1\over
  K^2}  \int_{|z|=1}\big|(F^\alpha(z))'\big|^2
  |{\rd}z|\\
  &\ll&
     {(e_rD(1))^{2\alpha}\over K^2}  \int_{|z|=1}|q(z)|^2
  |{\rd}z|\ll {(e_rD(1))^{2\alpha}n\over K^2}.
\end{eqnarray*}
  Collecting the last three estimates, by Lemma \ref{l-2}, we obtain
  the desired claim.

  Lemma \ref{l-5} is proved.
\bigskip

At this stage we have the following estimate.

\begin{lem} \label{l-6} If   Condition $(\ref{d})$ is satisfied and $\delta n\geq1$, then there exists a
positive constant $c= c(d')$ such that
\begin{equation}
    F_m=J_0+\rO\big(e_r D_n \delta^c\big)
\label{Fm}
\end{equation}
 uniformly in $0\leq r\leq \delta n$ and $n/2\leq m \leq
n$. Moreover,
\begin{equation}
      D_n={1\over 2\pi i n}\int_{\Delta_0} D'(z){\rd z\over z^n}+
      \rO\big(D_n\delta^c\big).
\label{Dn}
\end{equation}
\end{lem}

  {\it Proof}. It suffices to apply Lemmas \ref{l-3} and \ref{l-5} with
$K=\delta^{c(\alpha)}n $, where
\[
   c(\alpha)=\min\{1, \, d'(1-\alpha)/(\alpha d'+1)\},
   \]
and optimize the function $d'\alpha c(\alpha)$ with respect to
$\alpha\in (0,1)$. If $d'\leq 3$, then (\ref{Fm}) holds with
$c=(\sqrt{1+d'}-1)^2$.  If $d'>3$, the choice $\alpha=(d'-1)/2d'$
gives  $c(\alpha)=1$; thus, (\ref{Fm}) holds with $c=(d'-1)/2$. To
obtain (\ref{Dn}), use (\ref{Fm}) with $r=0$ and $m=n$.

  The lemma is proved.

\begin{lem} \label{l-7} If    $0\leq \eta\leq 1/2$ and $1/n\leq \delta\leq 1/2$ are arbitrary, then
\[
           J_0 = e_r D_n\bigg(1+\rO\Big(\big(\eta+(r/n){\bf1}\{r\geq 1\}\big) \delta^{-1}+\delta^c\Big)
           \bigg)
\]
uniformly in $n(1-\eta)\leq m \leq n$ and $0\leq r\leq \delta n$
with the constant $c$ defined in Lemma $\ref{l-5}$.
\end{lem}

{\it Proof}. If $z\in \Delta_0$ and $r\geq 1$, then
\begin{eqnarray*}
F'(z)&=& e_r D(z)\exp\bigg\{-\sum_{j\leq r}{d_j\over
j}(z^j-1)\bigg\}q(z)\\
&=& e_r D(z) \bigg(1+ \rO\Big({r\over\delta
n}\Big)\bigg)\bigg(\sum_{j\leq
n}-\sum_{j\leq r}\bigg)d_jz^{j-1} \\
&=& e_r D'(z) \Big(1+ \rO\big(r/\delta n\big)\Big)+\rO\big(r e_r
D(1)\big)
\end{eqnarray*}
and
\[
z^{-m}=z^{-n} \big(1+\rO(\eta \delta^{-1})\big).
\]
Consequently, by virtue of  $m^{-1}= n^{-1}\big(1+\rO(\eta)\big)$,
from Lemma \ref{l-2} and Equation (\ref{Dn}), we obtain
\begin{eqnarray*}
J_0&=&{e_r\over 2\pi i n}\bigg(1+ \rO\Big(\Big({r\over
n}+\eta\Big)\,{1\over \delta}\Big)\bigg)\int_{\Delta_0}
D'(z){\rd z\over z^n}+\rO\Big(e_rD_n {r\over \delta n}\Big)\\
&=&e_rD_n\Big(1+ \rO\big((r/n+\eta) \delta^{-1}+\delta^c\big)\Big).
\end{eqnarray*}

  If $r<1$, the terms having the fraction $r/n$ do not appear.

   The lemma is proved.

\smallskip

{\it Proof of Proposition \ref{prop1}}. Apply (\ref{Fm}) and the
last lemma.
\medskip

{\it Proof of Fundamental Lemma}. We now apply Lemma 1 and
Proposition~\ref{prop1} with $d_j=\lambda_j$. Condition (\ref{d})
for weakly logarithmic assemblies  is satisfied. From (\ref{P}) and
Proposition~\ref{prop1} with $\eta=(r/n)^{1/2}$ and
$\delta=(r/n)^{1/2(1+c)}$, we obtain
\[
{P\big(\ell_{rn}(\bar \xi)=n-m\big)\over
P\big(\ell(\bar\xi)=n\big)}=1+\rO\big((r/n)^{c_0}\big), \quad
c_0:=c/2(1+c),
\]
uniformly in $0\leq m\leq \sqrt{r n}$ provided that $1\leq r\leq
4^{-1-c} n$.

The summands over $m>\sqrt{r n}$ in (\ref{onedim}) contribute not
more than
\[
(r n)^{-1/2}\E\ell_r(\bar\xi)=(r n)^{-1/2}\sum_{j\leq r}
j\lambda_j\leq \theta'' (r/n)^{1/2}.
\]
Hence, by (\ref{onedim}), we obtain
\[
 \rho_{TV}\Big(\cL\big(\bar k_r(\sigma)\big),\cL(\bar
    \xi_r)\Big)\ll (r/n)^{c_1},
 \]
 where $c_1=\min\{1/2, c_0\}$ and $1\leq r\leq 4^{-1-c} n$. Since
 the claim of Fundamental Lemma is trivial for $n\leq 4^{1+c}$, we
 have finished its proof.

\section{Proof of Theorem 2 and its Corollaries}
\label{S-3} Set  $\Z_+^n(m) = \{\bar s\in \Z_+^n:\;\ell(\bar s)=m\}$
where $0\leq m\leq n$. For arbitrary distributions
 $p_j(k)$, $1\leq j\leq n$, on $\Z_+$ we define the product measure on
 $\Z_+^n$ by
 \[
             P(\{\bar k\})=\prod_{j\leq n} p_j(k_j), \quad \bar
             k=(k_1,\dots ,k_s)\in\Z_+^n.
 \]
Denote for brevity $P_n=P(\Z_+^n(n))$. Let $V=V(U)$ be the extension
of an arbitrary subset $U\subset\Z_+^n$ defined in (\ref{VU}).

\begin{lem} \label{l-8} Suppose $n\ge1$ and there exist positive
constants
         $c_2,c_3, C_1,C_2$ such that
\begin{itemize}
\item[$(i)$]  $p_j(0)\ge c_2$ for  all $1\le j\le n$ $;$

\item[$(ii)$] $P\big(\Z_+^n(m)\big)\le C_1\, \bigg(\displaystyle\frac
n{m+1}\bigg)^{1-\theta}
      P_n$ for  $0\le m\leq n-1$ and for some $0 < \theta \leq 1$
      $;$

\item[$(iii)$]  $ P_n\ge c_3n^{-1}$ $;$

\item[$(iv)$] for $1\le m\le n$,
    $$
         \displaystyle\sum_{\footnotesize\begin{array}{c} k\ge1,j\le n\\ kj=m\end{array}}
          \frac{p_j(k)}{p_j(0)}\le \frac{C_2}{m}.
    $$
\end{itemize}
Then
    $$
         P\big(\overline{V}\;|\Z_+^n(m)\big)\le CP^\theta(\overline{U})
         +C_1C_2\theta^{-1}n^{-\theta}{\bf 1}\{\theta<1\},
    $$
where
    $$
         C:=\max\bigg\{\frac{32}{c_2^2}, \frac{C_2}{c_3}+\frac{4C_1}{c_2}
         +\frac{C_1C_2}{\theta} \bigg\}.
    $$
\end{lem}

{\it Proof}.  See \cite{EJ-99}, Appendix.

\smallskip

{\it Proof of Theorem} 1. It suffices to check conditions $(i)-(iv)$
of the last lemma for the poissonian probabilities $p_j(k)$ with
parameters $\lambda_j$. By virtue of Condition (\ref{Def}), $(i)$
and $(iv)$ are trivial. Further, we find
\begin{eqnarray*}
P\big(\Z_+^n(m)\big)&=&P\bigg(\sum_{j=1}^m
j\xi_j=m,\xi_{m+1}=0, \dots,\xi_n=0\bigg)\nonumber\\
&=&\exp \bigg\{-\sum _{j=1}^{n}\lambda_j\bigg\}\sum_{\ell_m(\bar
k)=m}\prod_{j=1}^{m}{{\lambda_j^{k_j}}\over {k_j!}}\nonumber\\
&=&\exp \bigg\{-\sum _{j=1}^{n}\lambda_j\bigg\}
[z^m]\exp\bigg\{\sum_{j\leq m}\lambda_j z^j\bigg\}, \quad 0\leq
m\leq n.
\end{eqnarray*}
Hence,  applying Lemma \ref{l-2}, we obtain
\[
    P\big(\Z_+^n(m)\big)\asymp {1\over m+1}\exp \bigg\{-\sum
    _{j=m+1}^{n}\lambda_j\bigg\}
    \]
for $ 0\leq m\leq n$, where $a\asymp b$ means $a\ll b\ll a$. This
and Condition (\ref{Def}) imply $(ii)$ and $(iii)$.

The theorem is proved.

\s

{\it Proof of Corollary} 1. Apply Theorem \ref{T-2} for
\[
    U=\bigg\{\bar t\in\Z_+^n:\; H(\bar t)\in
    A\bigg\},
\]
where $H(\bar t):=\sum_{j\leq n} h_j(t_j)$, and check that
\[
   V(U)\subset\big\{ \bar s\in\Z_+^n:\; H(\bar s)\in A+A-A\big\}.
\]
Now \begin{eqnarray*}
   \mu_n\big(h(\sigma)\not\in A+A-A\big)&=&P\big(H(\bar\xi)\not\in A+A-A|\, \ell(\bar
   \xi)=n\big)\\
   &\leq& P\big(\bar\xi\not\in V(U)|\, \ell(\bar
   \xi)=n\big)\\
   &\ll& P^{\theta}\big(\bar\xi\not\in U\big)+{\bf
   1}\{\theta'<1\} n^{-\theta'}\\
   &=&P^{\theta}\big(H(\bar\xi)\not\in A\big)+{\bf
   1}\{\theta'<1\} n^{-\theta'}.
\end{eqnarray*}

  Corollary 1 is proved.
  \s

{\it Proof of Corollary} 2. Apply the previous corollary for $\G=\R$
and $A=\{t:\; |t-a|\leq u/3\}$.

\section{Proof of Theorem 2}

We adopt the argument used in the case of permutations
\cite{EM-98LMJ} and for the logarithmic assemblies \cite{JN1}.

Let $Z_1, Z_2,\dots, Z_n$ be independent random variables defined on
some probability space $(\Omega, {\cal F}, P)$, with
$\mathbb{E}Z_j=0, {\E}Z_j^2<\infty, j=1,2,\dots$, and
\[
     S_m=\sum_{j=1}^mZ_j, \qquad
     D^2_m=\sum_{j=1}^m {\E}Z_j^2.
\]
 We define the polygonal lines $s_n(\cdot):[0,D_n^2]\to \mathbb{R}$ such that
$$
s_n(t)=S_m{D^2_{m+1}-t\over D^2_{m+1}-D^2_m}+S_{m+1}{t-D^2_m\over
D^2_{m+1}-D^2_m}
$$
if $D_m^2\leq t<D_{m+1}^2$ and $ 0\leq m\leq n-1$. Set also
\[
       S_n(t)={s_n(D^2_nt)\over \sqrt{2D^2_nLLD^2_n}}
\]
 for $0\leq t\leq1$ and $n\in \mathbb{N}$.

\begin{lem} \label{l-Ma} Let $D(n)\to \infty$ as
$n\to \infty$. Assume that there exists a sequence
$$
M_n={\rm o}\bigg({D_n\over\sqrt{LL D_n}}\bigg)
$$
such that
$$
P\big(|Z_n|\le M_n\big)=1
$$
for each $n\geq1$. Then
$$
S_n(\cdot)\Rightarrow \cK\quad (P\textrm{-a.s.}).
$$
  \end{lem}

{\it Proof}. This is Major's Theorem \cite{Ma}. \s

We will apply Lemma \ref{l-Ma}  for
$Z_j=a_j\big(\eta_j-(1-\re^{-\lambda_j})\big)$, where
$\eta_j:={\bf1}\{\xi_j\geq1\}$ and $1\leq j\leq n$. Then
$D_n^2=B^2(n)$ and  Condition (\ref{c-Kol}) will be at our disposal.
To simplify the calculations, we introduce another sequence of
additive functions
$$
      \tilde{h}(\sigma, m):=\sum_{j=1}^ma_j{\bf1}\{k_j(\sigma)\ge
      1\},\quad   m\leq n.
$$
Let $\tilde{u}_m(\sigma,t)$ and $\widetilde{U}_m(\sigma, t) $ be the
combinatorial processes
 defined as $u_m(\sigma,t)$ and $U_m(\sigma, t)$ using $\tilde{h}(\sigma,m)$
 instead of $h(\sigma,m)$. Set also $Y_m=a_1\eta_1+\cdots+a_m\eta_m$ for $1\leq m\leq n$.

\begin{lem} \label{l-cad} For arbitrary $\varepsilon>0$,
 \begin{eqnarray}
\lim_{n_1\to \infty}\limsup_{n\to \infty}\mu_n\Big(\max_{n_1\leq
m\leq n}\rho\big(\widetilde{U}_m(\sigma, \cdot), U_m(\sigma, \cdot)
\big)\ge \varepsilon\Big)=0.
\end{eqnarray}
\end{lem}

{\it Proof}. If $j$ and $j'$ are the consecutive numbers from the
set $I:=\{j\leq m:\; a_j\not=0\}$, then, by virtue of the definition
of $u_m(\sigma,t)$,
\begin{eqnarray*}
&&\max\Big\{|\widetilde{U}_m(\sigma,t)-U_m(\sigma,t)|:{B^2(j)\over
B^2(m)}\le t\le
{B^2(j')\over B^2(m)}\Big\}\\
&\le &\beta^{-1}(m) \max\Big\{|\tilde{h}(\sigma, j)-h(\sigma,
j)|,|\tilde{h}(\sigma, j')-h(\sigma, j')|\Big\}.
\end{eqnarray*}
Hence
\begin{eqnarray*}
&&\mu_n\Big(\max_{n_1\leq m\leq n}\rho\big(\widetilde{U}_m(\sigma,
\cdot), U_m(\sigma, \cdot) \big)\geq \varepsilon\Big)\\
&\le& \mu_n\Big(\max_{n_1\leq m\le n}\max_{j\in I}|\tilde{h}(\sigma,
j)-h(\sigma,j)|\geq \varepsilon\beta(n_1)\Big)\\
&\le&
\mu_n\bigg(\sum_{j=1}^n\big|h_j(k_j(\sigma))-a_j\cdot{\bf1}\{k_j(\sigma)\geq
1\}\big|\geq \varepsilon \beta (n_1)\bigg)\\
&\ll& P^{\theta}\bigg(\sum_{j=1}^n\big|h_j(\xi_j)-a_j\eta_j\geq
1\}\big|\geq (\varepsilon/3) \beta (n_1)\bigg) +\ro(1).
\end{eqnarray*}
In the last step we applied  Corollary \ref{cor-2}.
 In its turn, if  $K>2$ is arbitrary, the probability appearing on the right-hand side  can
be majorized by
\begin{eqnarray*}
&&
  P\big(\exists
j\le K: \xi_j\geq K\big)+P\big(\exists j>K: \xi_j\geq
2\big)\\
&+&P\bigg( \sum_{j\leq
K}\big(|h_j(\xi_j)|+|a_j|\eta_j\big)\geq(\varepsilon/3) \beta(n_1),
\quad 2\leq \xi_j\leq K, \forall j\leq K\bigg).
\end{eqnarray*}
 Since $\beta(n_1)\to\infty$ as $n_1\to\infty$, the
last probability is negligible. The first two of them do not exceed
\begin{eqnarray*}
\sum_{j\leq K}\sum_{k\geq K}{e^{-\lambda_j}\lambda_j^k\over
k!}+\sum_{j\geq K}\sum_{k\geq 2}{e^{-\lambda_j}\lambda_j^k\over
k!}\ll K^{-1}.
\end{eqnarray*}

Collecting the estimates, since $K$ is arbitrary, we obtain the
desired claim of Lemma \ref{l-cad}.

\s

In the sequel, we use only the functions $\tilde h(\sigma, m)$ and
the processes $\widetilde{U}_m(\sigma,t)$ writing them without the
"tilde".

\begin{lem} \label{l-trunc}  Let $1\leq k\leq n$,  $0<b_n \leq b_{n-1}\leq\cdots \leq b_1$,
and $\varepsilon > 0$ be arbitrary. For $h=\tilde h$, if
$n\to\infty$, we have
\begin{eqnarray*}
&&\mu_n\Big(\max_{k\leq m \leq n}b_m\big|h(\sigma,m)-A(m)\big|\geq \varepsilon \Big)\\
&\ll&
 P^{\theta}\Big(\max_{k\leq m \leq n}b_m\big|Y_m-A(m)\big|\geq \varepsilon/3 \Big)+ \ro(1)\\
&\leq& 3^{2\theta}\e^{-2\theta} \bigg( b_k^2B^2(k)+\sum_{k\leq j\leq
n} b_j^2a_j^2\re^{-\lambda_j}(1-\re^{-\lambda_j})\bigg)^{\theta}+
\ro(1).
\end{eqnarray*}
\end{lem}

{\it Proof.} The first estimate follows from Corollary \ref{cor-1}
applied for $\G=\R^{n-r+1}$,
\[
    A=\big\{ (s_r,\dots, s_n)\in\R^{n-r+1}:\; \max_{r\leq
    m\leq n}|s_m-A(m)|<\e/3\big\},
\]
and
\[
    h(\sigma)=\big(h(\sigma,r),\dots, h(\sigma,n)\big).
\]

  The second inequality in Lemma \ref{l-trunc} is just a partial case of Theorem 13 in Chapter III of
  \cite{Pe}.

  The lemma is proved.
\s

Let $r$, $n_1\leq r\leq n$, be a parameter, $q:=\max\{j\in I:\; j\le
r\}$, and
$$
u_m^{(r)}(\sigma,t)=\left\{%
\begin{array}{ll}
    u_m(\sigma,t)\quad & \hbox{\rm if}\quad t\leq B^2(q), \\
    u_m(\sigma, B^2(q))\quad  & \hbox{\rm if}\quad t>B^2(q). \\
  \end{array}%
  \right.
  $$
Denote $U_m^{(r)}(\sigma, t):=u_m^{(r)}\big(\sigma,
B^2(m)t\big)/\beta(m)$. Similarly, let
$$
s_m^{(r)}(t)=\left\{%
\begin{array}{ll}
    s_m(t)\quad & \hbox{\rm if}\quad t\leq B^2(q), \\
    s_m\big(B^2(q)\big)\quad  & \hbox{\rm if}\quad t>B^2(q) \\
  \end{array}%
  \right.
  $$
 and $S_m^{(r)}(t)=
s_m^{(r)}\big(tB(m)\big)/\beta(m)$.

\begin{lem} \label{l-trunc1} There exists a sequence $r=r(n)$, $n_1\leq r=\ro(n)$, such that, for every $\e>0$,
\begin{equation}
\lim_{n_1\to\infty}\limsup_{n\to\infty}P\Big(\max_{n_1\leq m\leq
n}\rho\big(S_m(\cdot),S_m^{(r)}(\cdot)\big)
 \geq \varepsilon\Big)=0
\label{Pn1}
 \end{equation}
and
\[
\lim_{n_1\to\infty}\limsup_{n\to\infty}\mu_n\Big(\max_{n_1\leq m\leq
n}\rho\big(U_m(\sigma,\cdot),U_m^{(r)}(\sigma,\cdot)\big)
 \geq \varepsilon\Big)=0.
 \]
 \end{lem}

{\it Proof}. If $P_{n_1,n}(\e)$ denotes the probability in
(\ref{Pn1}) and $n_1\leq r\leq n$, then
\begin{eqnarray*}
&&P_{n_1,n}(\e)=  P\Big(\max_{r\leq m\leq
n}\rho\big(S_m(\cdot),S_m^{(r)}(\cdot)\big)
 \geq \varepsilon\Big)\\
 &=&
P\Big(\max_{r\leq m\leq n}{1\over \beta(m)}\sup\big\{\big|
s_m(t)-s_m(B^2(q))\big|:\; {B^2(q)\leq t\leq B^2(m)}\big\}
 \geq \varepsilon\Big)\\
&\leq& P\Big(\max_{r< m\leq
n}\beta^{-1}(m)\big|(Y_m-A(m))-(Y_r-A(r))\big|\geq
\varepsilon\Big)\\
&\leq& \e^{-2}\, {B^2(n)-B^2(r)\over \beta^2(r)}
\end{eqnarray*}
by the already mentioned  Theorem 13 \cite{Pe}, Chapter III.

The same argument and Lemma \ref{l-trunc} (applied in the case
$a_j\equiv 0$ if $j\leq r$) leads to the estimate
\begin{eqnarray*}
&&\mu_n\Big(\max_{n_1\leq m\leq
n}\rho\big(U_m(\sigma,\cdot),U_m^{(r)}(\sigma,\cdot)\big)
 \geq \varepsilon\Big)\\
 &\leq&
\mu_n\Big(\max_{r<m\leq
n}\beta^{-1}(m)\big|(h(\sigma,m)-A(q))-(h(\sigma,r)-A(r))\big|\geq
\varepsilon\Big)\\
&\ll& P^\theta\Big(\max_{r<m\leq
n}\beta^{-1}(m)\big|(Y_m-A(m))-(Y_r-A(r))\big|\geq
(1/3)\varepsilon\Big)+\ro(1)\\
&\ll&\bigg({B^2(n)-B^2(r)\over \beta^2(r)}\bigg)^\theta +\ro(1)
 \end{eqnarray*}
as $n\to\infty$.

By Condition (4), if $r$ is sufficiently large, $r\leq j\leq n$, and
 $\delta$, $0<\delta<1$, is arbitrary, then  $|a_j|\leq \delta
B(n)/\sqrt{LL B(n)}$. Hence,  taking $r=\delta n$ and applying
Condition (\ref{d}), we obtain
\[
B^2(n)-B^2(r)\ll\delta^2\log {1\over\delta}\,  {B^2(n)\over LL
B(n)}.
\]
We now choose $\delta=\delta_n=\ro(1)$ as $n\to\infty$ so that
$\delta\geq1/\sqrt n$. This implies $B^2(n)-B^2(r)=\ro(\beta^2(r))$.
Having in mind the above estimates, we see that, with such an  $r$,
the probabilities in Lemma \ref{l-trunc1} vanish as $n\to\infty$ and
$n_1\to\infty$.

  The lemma is proved.

\medskip

{\it Proof of Theorem} 2. By virtue of the definition of strong
convergence and Lemma \ref{l-trunc1}, it suffices to prove that
$$
\lim_{n_1\to \infty}\limsup_{n\to \infty}\mu_n\Big(\max_{n_1\le m\le
n}\rho\big(U_m^r(\sigma,\cdot), \cK\big)\ge \varepsilon \Big)=0
$$
and
$$
\lim_{n_1\to \infty}\liminf_{n\to\infty}\mu_n\Big(\min_{n_1\le m\le
n}\rho\big(U_m^r(\sigma,\cdot),  g\big)<\e\Big)=1
$$
for each  function $g\in\cK$ and $\e>0$. Since here
$r=r(n)\to\infty$ and $r=\ro(n)$, we can apply the Fundamental Lemma
and substitute the frequencies by the appropriate probabilities for
independent r.vs. Consequently, our task reduces to the proof of
\[
\lim_{n_1\to \infty}\limsup_{n\to \infty}P\Big(\max_{n_1\le m\le
n}\rho\big(S_m^r(\cdot), \cK\big)\ge \varepsilon \Big)=0
\]
and
\[
\lim_{n_1\to \infty}\liminf_{n\to\infty}P\Big(\min_{n_1\le m\le
n}\rho\big(S_m^r(\cdot),  g\big)<\e\Big) =1.
\]
Checking that the last relations follow from  Lemmas \ref{l-Ma} and
\ref{l-trunc1} we complete the proof of Theorem 2.




\begin{thebibliography}{99}

\bibitem{ABT1}  R.\ Arratia, A.~D.~ Barbour and S.\ Tavar\'e, Limits of  logarithmic combinatorial
     structures, Ann. Probab. 28 (2000) 1620--1644.

\bibitem{ABT} R.\ Arratia, A.\ D.\ Barbour and S.\ Tavar\'e, Logarithmic
        Combinatorial Structures: a Probabilistic Approach, EMS
        Monographs in Mathematics, EMS Publishing House, Z\"urich, 2003.

\bibitem{EJ-99} G.\ J.\ Babu and E.\ Manstavi\v cius,  Brownian motion and
        random permutations,  Sankhy$\bar{a}$ A 61 (1999) 312--327.

\bibitem{EJ-02} G.\ J.\ Babu and E.\ Manstavi\v cius, Limit processes with independent
increments for the Ewens sampling formula, Ann. Inst. Statist. Math.
54 (2002) 607--620.

 \bibitem{EJ-02a} G.\ J.\ Babu and E.\ Manstavi\v cius, Infinitely divisible limit processes for the Ewens sampling
formula,  Lith. Math. J. 42 (2002) 232--242.


\bibitem{J-E-V-07} G.\ J.\ Babu, E.\ Manstavi\v cius and V. Zacharovas,
Limiting processes with dependent increments for measures on
symmetric group of permutations, in: S. Akiyama  {\it et al} (Eds)
Probability and Number Theory - Kanazawa 2005,  Advanced Studies in
Pure Math. 49, Math. Soc. Japan, Tokyo, 2007, pp. 41--67.

\bibitem{Fel} W.~Feller, The general form of the so-called law of
the iterated logarithm,  Trans. Amer. Math. Soc. 54 (1943) 373--402.


\bibitem{F-S} P. \ Flajolet and M. \ Soria, Gaussian limiting
distributions for the number of components in combinatorial
     structures,  J. Combin. Theory Ser. A  53 (1990) 165--182.

\bibitem{F} D.~Foata,  La s\'erie g\'en\'eratrice exponentielle
dans les probl\`emes d'\'enum\'eration, S\'eminaire de
Math\'ematiques Sup\'erieures. Les Presses de l'Universit\'e de
Montr\'eal, Qu\'ebec, 1974.

\bibitem{Fr} D.~Freedman,  Brownian Motion and Diffusion, Holden--day San
Francisco, 1971.

\bibitem{H} J.\ Hansen, Order statistics for decomposable combinatorial
        structures,  Random Structures Algorithms  5
        (1994) 517--533.

\bibitem{Ma} P.~Major, A note on Kolmogorov's law of iterated logarithm.
 Stud. Scient. Math. Hung. 12 (1977) 191--167.


\bibitem{EM-96} E.\ Manstavi\v cius, Additive and multiplicative
        functions on random  permutations,   Lith. Math. J.
        36 (1996) 400--408.

\bibitem{EM-98} E.\ Manstavi\v cius, The Berry--Esseen bound in the theory of
        random permutations,  Ramanujan J. 2 (1998)
        185--199.

\bibitem{EM-98LMJ} E.\ Manstavi\v cius, The law of iterated logarithm for
           random permutations, Lith. Math. J. 38
           (1998) 160--171.

\bibitem{EM-02Pal} E.\ Manstavi\v cius,  Functional limit theorems for sequences of mappings on
the symmetric group, in: A.~Dubickas {\it et al} (Eds), Analytic and
Probab. Methods  in Number Theory,  TEV, Vilnius, 2002, pp. 175-187.

\bibitem{EM-02} E.\ Manstavi\v cius,  Mappings on decomposable combinatorial
        structures: analytic approach,  Combinatorics, Probab.
        Computing 11 (2002) 61--78.


\bibitem{EM-04}  E.~Manstavi\v cius, Iterated logarithm laws and the cycle lengths
of a random permutation, in: M.~Drmota {\it et al} (Eds), Trends
Math., Mathematics and Computer Science III, Algorithms, Trees,
Combinatorics and Probabilities,  Birkh\"auser, Basel, 2004, pp.
39--47.

\bibitem{EM-08} E.\ Manstavi\v cius,  Asymptotic value distribution of additive
         function defined on the symmetric group,  Ramanujan J.  17 (2008) 259--280.

\bibitem{EM-09} E.\ Manstavi\v cius, An analytic method in probabilistic
combinatorics,  Osaka J. Math. 46 (2009), 273--290.

\bibitem{EM-JN} E.\ Manstavi\v cius and J.\ Nork\=unien\.e,
An analogue of Feller's theorem for logarithmic combinatorial
assemblies,  Lith. Math. J.  48 (2008) 405--417.

\bibitem{JN} J.~Nork\=unien\.e, The law of iterated logarithm for
combinatorial assemblies, Lith. Math. J.  46 (2006)  432--445.

\bibitem{JN1} J.~Nork\=unien\.e, The Strassen law of iterated logarithm for combinatorial
assemblies,  Lith. Math. J. 47 (2007) 176--183.

\bibitem{Pe} V.V.~Petrov,  Sums of Independent Random
Variables,  Ergebnisse der Mathematik und ihrer Grenzgebiete, B. 82,
Springer, New York, 1975.


\bibitem{Ru} I.Z.~Ruzsa, Generalized moments of additive functions,
J. Number Theory 18 (1984) 27--33.

\bibitem{St1} D.~Stark, Explicit limits of total variation distance
in approximations of random logarithmic assemblies by related
Poisson processes,  Combinatorics, Probab.
        Computing  6 (1997) 87--105.

   \bibitem{St2} D.~Stark, Total variation asymptotics for related
Poisson process approximations of random logarithmic assemblies,
Combinatorics, Probab.
        Computing  8 (1999) 567--598.

   \bibitem{VZ1} V.~Zacharovas, The convergence rate in CLT for random variables on
permutations, in: A.~Dubickas {\it et al} (Eds), Analytic and
Probab. Methods  in Number Theory,
 TEV, Vilnius, 2002, pp. 329-338.

\bibitem{VZ2} V.~Zacharovas,  The convergence rate to the normal law of a certain
variable defined on random polynomials,  Lith. Math. J. 42 (2002)
88--107.

\bibitem{VZ3} V.~Zacharovas, Distribution of the logarithm of the order of a random
permutation,  Lith. Math. J. 44 (2004) 296--327.

\bibitem{VZ4}  V.~Zacharovas, Distribution of random variables on the symmetric
group, Doctorial dissertation, Vilnius University, 2004;
arXiv:0901.1733.
\end{thebibliography}
\end{document}